\documentclass[fleqn]{llncs}
\usepackage{mepa}
\usepackage{version}

\title{Meadow Enriched ACP Process Algebras}
\author{J.A. Bergstra \and C.A. Middelburg}
\institute{Programming Research Group, University of Amsterdam, \\
           P.O.~Box~41882, 1009~DB~Amsterdam, the Netherlands \\
           \email{J.A.Bergstra@uva.nl,C.A.Middelburg@uva.nl}}

\begin{document}

\maketitle

\begin{abstract}
We introduce the notion of an \ACP\ process algebra.
The models of the axiom system \ACP\ are the origin of this notion.
\ACP\ process algebras have to do with processes in which no data are
involved.
We also introduce the notion of a meadow enriched \ACP\ process algebra,
which is a simple generalization of the notion of an \ACP\ process
algebra to processes in which data are involved.
In meadow enriched \ACP\ process algebras, the mathematical structure
for data is a meadow.
\begin{classcode}
D.1.3, D.2.1, D.2.4, F.1.2, F.3.1.
\end{classcode}
\end{abstract}

\section{Introduction}
\label{sect-introduction}

In many formalisms proposed for the description and analysis of
processes in which data are involved, algebraic specifications of
the data types concerned have to be given over and over again.
This is also the case with the principal \ACP-based formalisms proposed
for the description and analysis of processes in which data are
involved, to wit $\mu$CRL~\cite{GP94a,GP95a} and PSF~\cite{MV90a}.
There is a mismatch between the process specification part and the data
specification part of these formalisms.
Firstly, there is a choice of one built-in type of processes, whereas
there is a choice of all types of data that can be specified
algebraically.
Secondly, the semantics of the data specification part is its initial
algebra in the case of PSF and its class of minimal Boolean preserving
algebras in the case of $\mu$CRL, whereas the semantics of the process
specification part is a model based on transition systems and
bisimulation equivalence.
Sticking to this mismatch, no stable axiomatizations in the style of
\ACP\ has emerged for process algebras that have to do with processes
in which data is involved.

Our objective is to obtain a stable axiomatization in the style of \ACP\
for process algebras that have to do with processes in which data is
involved.
To achieve this objective, we first introduce the notion of an \ACP\
process algebra and then the notion of a meadow enriched \ACP\ process
algebra.
\ACP\ process algebras are essentially models of the axiom system \ACP.
Meadow enriched \ACP\ process algebras are data enriched \ACP\ process
algebras in which the mathematical structure for data is a meadow.
Meadows are defined for the first time in~\cite{BT07a}.
The prime example of meadows is the rational number field with the
multiplicative inverse operation made total by imposing that the
multiplicative inverse of zero is zero.
Although the notion of a meadow enriched \ACP\ process algebra is a
simple generalization of the notion of an \ACP\ process algebra, it is
an interesting one: there is a multitude of finite and infinite meadows
and meadows obviate the need for Boolean values and operations on data
that yield Boolean values to deal with conditions on data.

In the work on \ACP, the emphasis has always been on axiom systems.
In this paper, we put the emphasis on algebras.
That is, \ACP\ process algebras are looked upon in the same way as
groups, rings, fields, etc.\ are looked upon in universal algebra (see
e.g.~\cite{BS81a}).
The set of equations that are taken to characterize \ACP\ process
algebras is a revision of the axiom system \ACP.
The revision is primarily a matter of streamlining.
However, it also involves a minor generalization that allows for the
generalization to meadow enriched \ACP\ process algebras to proceed
smoothly.

This paper is organized as follows.
First, we give a brief summary of meadows (Section~\ref{sect-meadows}).
Next, we introduce the notion of an \ACP\ process algebra
(Section~\ref{sect-ACP-process-alg}) and the notion of an meadow
enriched \ACP\ process algebra (Section~\ref{sect-MD-ACP-process-alg}).
Finally, we make some concluding remarks
(Section~\ref{sect-conclusions}).

\section{Meadows}
\label{sect-meadows}

In the data enriched \ACP\ process algebras introduced in this paper,
the mathematical structure for data is a meadow.
A meadow is a field with the multiplicative inverse operation made total
by imposing that the multiplicative inverse of zero is zero.
Meadows are defined for the first time in~\cite{BT07a} and are
investigated in e.g.~\cite{BHT09a,BP08a,BR07a}.
In this section, we give a brief summary of meadows.

The signature of meadows is the same as the signature of fields.
It is a one-sorted signature.
We make the single sort explicit because we will extend this signature
to a two-sorted signature in Section~\ref{sect-MD-ACP-process-alg}.
The signature of meadows consists of the sort $\Quant$ of
\emph{quantities} and the following constants and operators:
\begin{itemize}
\item
the constants $\const{0}{\Quant}$ and $\const{1}{\Quant}$;
\item
the binary \emph{addition} operator
$\funct{+}{\Quant \x \Quant}{\Quant}$;
\item
the binary \emph{multiplication} operator
$\funct{\mul}{\Quant \x \Quant}{\Quant}$;
\item
the unary \emph{additive inverse} operator
$\funct{-}{\Quant}{\Quant}$;
\item
the unary \emph{multiplicative inverse} operator
$\funct{\minv}{\Quant}{\Quant}$.
\end{itemize}

We assume that there is a countably infinite set $\cU$ of variables,
which contains $u$, $v$ and $w$, with and without subscripts.
Terms are build as usual.
We use infix notation for the binary operators ${} + {}$ and
${} \mul {}$, prefix notation for the unary operator ${} -$, and postfix
notation for the unary operator ${}\minv$.
We use the usual precedence convention to reduce the need for
parentheses.
We introduce subtraction and division as abbreviations:
$p - q$ abbreviates $p + (-q)$ and
$p / q$ abbreviates $p \mul q\minv$.
We freely use the numerals $2,3,\ldots$\ .

The constants and operators from the signature of meadows are adopted
from rational arithmetic, which gives an appropriate intuition about
these constants and operators.

A meadow is an algebra with the signature of meadows that satisfies
the equations given in Table~\ref{eqns-meadow}.%
\begin{table}[!t]
\caption{Axioms for meadows}
\label{eqns-meadow}
\begin{eqntbl}
\begin{eqncol}
(u + v) + w = u + (v + w)                                             \\
u + v = v + u                                                         \\
u + 0 = u                                                             \\
u + (-u) = 0
\end{eqncol}
\qquad\quad
\begin{eqncol}
(u \mul v) \mul w = u \mul (v \mul w)                                 \\
u \mul v = v \mul u                                                   \\
u \mul 1 = u                                                          \\
u \mul (v + w) = u \mul v + u \mul w
\end{eqncol}
\qquad\quad
\begin{eqncol}
{(u\minv)}\minv = u                                                   \\
u \mul (u \mul u\minv) = u
\end{eqncol}
\end{eqntbl}
\end{table}
Thus, a meadow is a commutative ring with identity equipped with a
multiplicative inverse operation ${}\minv$ satisfying the reflexivity
equation ${(u\minv)}\minv = u $ and the restricted inverse equation
$u \mul (u \mul u\minv) = u$.
From the equations given in Table~\ref{eqns-meadow}, the equation
$0\minv = 0$ is derivable.

A \emph{non-trivial meadow} is a meadow that satisfies
the \emph{separation axiom}
\begin{ldispl}
0 \neq 1\;.
\end{ldispl}
A \emph{cancellation meadow} is a meadow that satisfies the
\emph{cancellation axiom}
\begin{ldispl}
u \neq 0 \And u \mul v = u \mul w \Implies v = w\;,
\end{ldispl}
or equivalently, the \emph{general inverse law}
\begin{ldispl}
u \neq 0 \Implies u \mul u\minv = 1\;.
\end{ldispl}
In~\cite{BT07a}, cancellation meadows are called zero-totalized fields.
An important  property of cancellation meadows is the following:
$0 / 0 = 0$, whereas $u / u = 1$ for $u \neq 0$.

\section{\ACP\ Process Algebras}
\label{sect-ACP-process-alg}

In this section, we introduce the notion of an \ACP\ process algebra.
This notion originates from the models of the axiom system \ACP, which
was first presented in~\cite{BK84b}.
A comprehensive introduction to \ACP\ can be found in~\cite{BW90,Fok00}.

It is assumed that a fixed but arbitrary finite set $\Act$ of
\emph{atomic action names}, with $\dead \notin \Act$, has been given.

The signature of \ACP\ process algebras is a one-sorted signature.
We make the single sort explicit because we will extend this signature
to a two-sorted signature in Section~\ref{sect-MD-ACP-process-alg}.
The signature of \ACP\ process algebras consists of the sort $\Proc$ of
\emph{processes} and the following constants, operators, and predicate
symbols:
\begin{itemize}
\item
the \emph{deadlock} constant $\const{\dead}{\Proc}$;
\item
for each $e \in \Act$,
the \emph{atomic action} constant $\const{e}{\Proc}$;
\item
the binary \emph{alternative composition} operator
$\funct{\altc}{\Proc \x \Proc}{\Proc}$;
\item
the binary \emph{sequential composition} operator
$\funct{\seqc}{\Proc \x \Proc}{\Proc}$;
\item
the binary \emph{parallel composition} operator
$\funct{\parc}{\Proc \x \Proc}{\Proc}$;
\item
the binary \emph{left merge} operator
$\funct{\leftm}{\Proc \x \Proc}{\Proc}$;
\item
the binary \emph{communication merge} operator
$\funct{\commm}{\Proc \x \Proc}{\Proc}$;
\item
for each $H \subseteq \Act$,
the unary \emph{encapsulation} operator
$\funct{\encap{H}}{\Proc}{\Proc}$;
\item
the unary \emph{atomic action} predicate symbol
$\predt{\isact}{\Proc}$.
\end{itemize}

We assume that there is a countably infinite set $\cX$ of variables of
sort $\Proc$, which contains $x$, $y$ and $z$, with and without
subscripts.
Terms are built as usual.
We use infix notation for the binary operators.
We use the following precedence conventions to reduce the need for
parentheses: the operator $\altc$ binds weaker than all other binary
operators and the operator $\seqc$ binds stronger than all other binary
operators.
We write
$\Altc{i \in \cI} P_i$, where $\cI = \set{i_1,\ldots,i_n}$ and
$P_{i_1},\ldots,P_{i_n}$ are terms of sort $\Proc$, for
$P_{i_1} \altc \ldots \altc P_{i_n}$.
The convention is that $\Altc{i \in \cI} P_i$ stands for $\dead$ if
$\cI = \emptyset$.

Let $P$ and $Q$ be closed terms of sort $\Proc$.
Intuitively, the constants and operators introduced above
can be explained as follows:
\begin{itemize}
\item
$\dead$ can neither perform an atomic action nor terminate successfully;
\item
$e$ first performs atomic action $e$ and then terminates successfully;
\item
$P \altc Q$ behaves either as $P$ or as $Q$, but not both;
\item
$P \seqc Q$ first behaves as $P$ and on successful termination of $P$ it
next behaves as $Q$;
\item
$P \parc Q$ behaves as the process that proceeds with $P$ and $Q$ in
parallel;
\item
$P \leftm Q$ behaves the same as $P \parc Q$, except that it starts
with performing an atomic action of $P$;
\item
$P \commm Q$ behaves the same as $P \parc Q$, except that it starts with
performing an atomic action of $P$ and an atomic action of $Q$
synchronously;
\item
$\encap{H}(P)$ behaves the same as $P$, except that atomic actions from
$H$ are blocked;
\item
$\isact(P)$ holds if $P$ is an atomic action.
\end{itemize}

The predicate $\isact$ is a means to distinguish atomic actions from
other processes.
An alternate way to accomplish this is to have a subsort $\mathbf{A}$ of
the sort $\Proc$.
We have not chosen this alternate way because it complicates matters
considerably.
Moreover, we prefer to keep close to elementary algebraic specification
(see e.g.~\cite{BT06a}).
By the notational convention introduced below, we seldom have to use
the predicate $\isact$ explicitly.

In equations between terms of sort $\Proc$, we will use a notational
convention which requires the following assumption: there is a countably
infinite set $\cX' \subseteq \cX$ that contains $a$, $b$ and $c$, with
and without subscripts, but does not contain $x$, $y$ and $z$, with and
without subscripts.
Let $\phi$ be an equation between terms of sort $\Proc$, and let
$\set{a_1,\ldots,a_n}$ be the set of all variables from $\cX'$ that
occur in $\phi$.
Then we write $\phi$ for
$\isact(x_1) \And \ldots \And \isact(x_n) \Implies \phi'$, where
$\phi'$ is $\phi$ with, for all $i \in [1,n]$, all occurrences of $a_i$
replaced by $x_i$, and
$x_1,\ldots,x_n$ are variables from $\cX$ that do not occur in $\phi$.

An \ACP\ process algebra is an algebra with the signature of \ACP\
process algebras that satisfies the formulas given in
Table~\ref{eqns-ACP}.%
\begin{table}[!t]
\caption{Axioms for \ACP\ process algebras}
\label{eqns-ACP}
\begin{eqntbl}
\begin{eqncol}
x \altc y = y \altc x                                           \\
(x \altc y) \altc z = x \altc (y \altc z)                       \\
x \altc x = x                                                   \\
(x \altc y) \seqc z = x \seqc z \altc y \seqc z                 \\
(x \seqc y) \seqc z = x \seqc (y \seqc z)                       \\
x \altc \dead = x                                               \\
\dead \seqc x = \dead                                           \\
{}                                                              \\
\encap{H}(e) = e     \hfill \mif e \notin H                     \\
\encap{H}(e) = \dead \hfill \mif e \in H                        \\
\encap{H}(\dead) = \dead                                        \\
\encap{H}(x \altc y) = \encap{H}(x) \altc \encap{H}(y)          \\
\encap{H}(x \seqc y) = \encap{H}(x) \seqc \encap{H}(y)
\end{eqncol}
\qquad\quad
\begin{eqncol}
x \parc y =
        (x \leftm y \altc y \leftm x) \altc x \commm y          \\
a \leftm x = a \seqc x                                          \\
a \seqc x \leftm y = a \seqc (x \parc y)                        \\
(x \altc y) \leftm z = x \leftm z \altc y \leftm z              \\
a \commm b \seqc x = (a \commm b) \seqc x                       \\
a \seqc x \commm b \seqc y = (a \commm b) \seqc (x \parc y)     \\
(x \altc y) \commm z = x \commm z \altc y \commm z              \\
x \commm y = y \commm x                                         \\
(x \commm y) \commm z = x \commm (y \commm z)                   \\
\dead \commm x = \dead                                          \\
{}                                                              \\
\isact(e)                                                       \\
\isact(x) \And \isact(y) \Implies \isact(x \commm y)
\end{eqncol}
\end{eqntbl}
\end{table}
Three equations in this table are actually schemas of equations:
$e$ is a syntactic variable which stands for an arbitrary constant of
sort $\Proc$ different from $\dead$.
A side condition is added to two schemas to restrict the constants for
which the syntactic variable stands.
The number of proper equations is still finite because there exists only
a finite number of constants of sort $\Proc$.

Because the notational convention introduced above is used, the four
equations in Table~\ref{eqns-ACP} that are actually conditional
equations look the same as their counterpart in the axiom system \ACP.
It happens that these conditional equations allow for the generalization
to meadow enriched \ACP\ process algebras to proceed smoothly.
Apart from this, the set of formulas given in Table~\ref{eqns-ACP}
differs from the axiom system \ACP\ on three points.
Firstly, the equations $x \commm y = y \commm x$,
$(x \commm y) \commm z = x \commm (y \commm z)$, and
$\dead \commm x = \dead$ have been added.
In the axiom system \ACP, all closed substitution instances of these
equations are derivable.
Secondly, the equations $a \seqc x \commm b = (a \commm b) \seqc x$
and $x \commm (y \altc z) = x \commm y \altc x \commm z$ have been
removed.
These equations can be derived using the added equation
$x \commm y = y \commm x$.
Thirdly,\linebreak[2] the formulas $\isact(e)$ and
$\isact(x) \And \isact(y) \Implies \isact(x \commm y)$ have been added.
They express that the processes denoted by constants of sort $\Proc$ are
atomic actions and that the processes that result from the communication
merge of two atomic actions are atomic actions.
This does not exclude that there are additional atomic actions, which is
impossible in the case of \ACP.

Not all processes in a \ACP\ process algebra have to be interpretations
of closed terms, even if all atomic actions are interpretations of
closed terms.
The processes concerned may be solutions of sets of recursion equations.
It is customary to restrict the attention to \ACP\ process algebras
satisfying additional axioms by which sets of recursion equations that
fulfil a guardedness condition have unique solutions.
For an comprehensive treatment of this issue, the reader is referred
to~\cite{BW90}.

\section{Meadow Enriched \ACP\ Process Algebras}
\label{sect-MD-ACP-process-alg}

In this section, we introduce the notion of an meadow enriched \ACP\
process algebra.
This notion is a simple generalization of the notion of an \ACP\
process algebra introduced in Section~\ref{sect-ACP-process-alg} to
processes in which data are involved.
The elements of a meadow are taken as data.

The signature of meadow enriched \ACP\ process algebras is a two-sorted
signature.
It consists of the sorts, constants and operators from the signatures of
\ACP\ process algebras and meadows and in addition the following
operators:
\begin{itemize}
\item
for each $n \in \Nat$ and $e \in \Act$,
the $n$-ary \emph{data handling atomic action} operator
$\funct{e}
 {\underbrace{\Quant \x \cdots \x \Quant}_{n\; \mathrm{times}}}{\Proc}$;
\item
the binary \emph{guarded command} operator
$\funct{\gc}{\Quant \x \Proc}{\Proc}$.
\end{itemize}

We take the variables in $\cU$ for the variables of sort $\Quant$.
We assume that the sets $\cU$ and $\cX$ are disjunct.
Terms are built as usual for a many-sorted signature
(see e.g.~\cite{ST99a,Wir90a}).
We use the same notational conventions as before.
In addition, we use infix notation for the binary operator ${} \gc {}$.

Let $p_1,\ldots,p_n$ and $p$ be closed terms of sort $\Quant$ and $P$ be
a closed term of sort $\Proc$.
Intuitively, the additional operators can be explained as follows:
\begin{itemize}
\item
$e(p_1,\ldots,p_n)$ first performs data handling atomic action
$e(p_1,\ldots,p_n)$ and then terminates successfully;
\item
$p \gc P$ behaves as the process $P$ under the condition that
$p = 0$ holds.
\end{itemize}

The different guarded command operators that have been proposed before
in the setting of \ACP\ have one thing in common: their first operand is
considered to stand for an element of the domain of a Boolean algebra
(see e.g.~\cite{BM05a}).
In contrast with those guarded command operators, the first operand of
the guarded command operator introduced here is considered to stand for
an element of the domain of a meadow.

A meadow enriched \ACP\ process algebra is an algebra with the signature
of meadow enriched \ACP\ process algebras that satisfies the formulas
given in Tables~\ref{eqns-ACP} and~\ref{eqns-ACPmd}.%
\begin{table}[!t]
\caption{Additional axioms for meadow enriched \ACP\ process algebras}
\label{eqns-ACPmd}
\begin{eqntbl}
\begin{array}{@{}l@{}}
\begin{eqncol}
0 \gc x = x                                                       \\
1 \gc x = \dead                                                   \\
u \gc x = u / u \gc x                                             \\
u \gc (v \gc x) = (1 - (1 - u/u) \mul (1 - v/v)) \gc x            \\
u \gc x \altc v \gc x = (u/u \mul v/v) \gc x                      \\
\end{eqncol}
\qquad\quad
\begin{eqncol}
u \gc \dead = \dead                                               \\
u \gc (x \altc y) = u \gc x \altc u \gc y                         \\
u \gc x \seqc y = (u \gc x) \seqc y                               \\
(u \gc x) \leftm y = u \gc (x \leftm y)                           \\
(u \gc x) \commm y = u \gc (x \commm y)                           \\
\encap{H}(u \gc x) = u \gc \encap{H}(x)
\end{eqncol}
\\
\begin{ceqncol}
e \commm e' = e'' \Implies
{} \\ \;\;
e(u_1,\ldots,u_n) \commm e'(v_1,\ldots,v_n) =
(u_1 - v_1) \gc
(\cdots \gc ((u_n - v_n) \gc e''(u_1,\ldots,u_n))\cdots)          \\
e \commm e' = \dead \Implies
e(u_1,\ldots,u_n) \commm e'(v_1,\ldots,v_n) = \dead               \\
e(u_1,\ldots,u_n) \commm e'(v_1,\ldots,v_m) = \dead
 & \hsp{-23} \mif n \neq m
\eqnsep
\encap{H}(e(u_1,\ldots,u_n)) = e(u_1,\ldots,u_n)
 & \hsp{-23} \mif e \not\in H                                   \\
\encap{H}(e(u_1,\ldots,u_n)) = \dead
 & \hsp{-23} \mif e \in H
\eqnsep
\isact(e(u_1,\ldots,u_n))
\end{ceqncol}
\end{array}
\end{eqntbl}
\end{table}
Like in Table~\ref{eqns-ACP}, some formulas in Table~\ref{eqns-ACPmd}
are actually schemas of formulas: $e$, $e'$ and $e''$ are syntactic
variables which stand for arbitrary constants of sort $\Proc$ different
from $\dead$ and, in addition, $n$ and $m$ stand for arbitrary natural
numbers.

For meadow enriched \ACP\ process algebras that satisfy the cancellation
axiom, the first five equations concerning the guarded command operator
can easily be understood by taking the view that $0$ and $1$ represent
the Boolean values $\True$ and $\False$, respectively.
In that case, we have that
\begin{itemize}
\item
$p/p$ models the test that yields $\True$ if $p = 0$ and $\False$
otherwise;
\item
if both $p$ and $q$ are equal to $0$ or $1$, then
$1 - p$ models $\Not p$, $p \mul q$ models $p \Or q$, and
consequently $1 - (1 - p) \mul (1 - q)$ models $p \And q$.
\end{itemize}
From this view, the equations given in the upper half of
Table~\ref{eqns-ACPmd} differ from the axioms for the most
general kind of guarded command operator that has been proposed in the
setting of \ACP\ (see e.g.~\cite{BM05a}) on two points only.
Firstly, the equation $u \gc x = u / u \gc x$ has been added.
This equation formalizes the informal explanation of the guarded command
given above.
Secondly, the equation $x \commm (u \gc y) = u \gc (x \commm y)$ has
been removed.
This equation can be derived using the equation
$x \commm y = y \commm x$ from Table~\ref{eqns-ACP}.

The equations in Table~\ref{eqns-ACPmd} concerning the communication
merge of data handling atomic actions formalize the intuition that two
data handling atomic actions $e(p_1,\ldots,p_n)$ and
$e'(q_1,\ldots,q_m)$ cannot be performed synchronously if $n \neq m$ or
$p_1 \neq q_1$ or $\ldots$ or $p_n \neq q_n$.
The equations concerning the encapsulation of data handling atomic
actions formalize the way in which it is dealt with in $\mu$CRL and PSF.
The formula concerning the atomic action predicate simply expresses that
data handling atomic actions are also atomic actions.

\section{Conclusions}
\label{sect-conclusions}

We have introduced the notion of an \ACP\ process algebra.
The set of equations that have been taken to characterize \ACP\ process
algebras is a revision of the axiom system \ACP.
We consider this revision worth mentioning of itself.
We have also introduced the notion of a meadow enriched \ACP\ process
algebra.
This notion is a simple generalization of the notion of an \ACP\ process
algebra to processes in which data are involved, the mathematical
structure of data being a meadow.
The primary mathematical structure for calculations is unquestionably a
field, and a meadow differs from a field only in that the multiplicative
inverse operation is made total by imposing that the multiplicative
inverse of zero is zero.
Therefore, we consider a combination of \ACP\ process algebras and
meadows like the one made in this paper, a combination with potentially
many applications.

\bibliographystyle{spmpsci}
\bibliography{PA}


\end{document}